\documentstyle[12pt]{article}
\font\erlss=eufm9
\newcommand{\Cs}{\mbox{\erlss C}}
\font\erl=eufm10 at 12pt
\newcommand{\C}{\mbox{\erl C}}
\begin{document}
%
%\baselineskip=7mm
%
%\vspace*{15mm}
%
%
\begin{flushleft}
{\large \bf Expansions of ${}_4\!\!\mbox{\boldmath$F$}_3$ When
the Upper
Parameters Differ by Integers}\\[5mm]
{\sc Megumi Saigo}\\
{\it Department of Applied Mathematics$,$ Fukuoka University$,$ 
Fukuoka 814-0180$,$ Japan\\[3mm]
{\sc R.K. Saxena}\\
Department of Mathematics and Statistics$,$ Jai Narain Vyas University$,$ 
Jodhpur-342001$,$ India}\\[3mm]
(1991 AMS Classification number: 33C20)\\[8mm]
\end{flushleft}
%
%\vspace{6mm}
%
{\small In this article three expansion formulas for a
generalized 
hypergeometric function $_4F_3$ are derived, when its upper
parameters 
differ by integers. Though the results are special cases of a 
general continuation formula for $_pF_q$, they are sufficiently
general 
and unify a number of known results.}\\[5mm]
\begin{flushleft}
{\bf 1. Introduction}
\end{flushleft}
\setcounter{section}{1}
\par
In many problems arising in physical sciences and statistics,
hypergeometric 
functions, in logarithmic cases, are applicable, which is known
from the 
monograph of Mathai and Haubold [4] and Mathai [3], etc. A
detailed account 
of such applications is available from the two monographs by
Mathai and 
Saxena [5], [6]. Expansions of the Gauss's hypergeometric
functions in 
logarithmic cases are given in the monograph of Erd\'elyi {\it et
al}. [1].
\par
In a series of papers by Saigo {\it et al}. [9] - [25], numerous
properties 
have been investigated which exhibits the behaviour of
generalized 
hypergeometric functions near the boundaries of their regions of
convergence 
of series defining these functions. In deriving some of these
results, the 
expansions of Gauss's hypergeometric function in logarithmic
cases given in 
the monograph by Magnus {\it et al.} [2] are used.
\par
The object of this article is to derive three expansions of a
generalized 
hypergeometric function $_4F_3(\cdot)$, when the upper parameters
differ by 
integers. Though the results established in this article are
special cases of a 
general continuation formula for $_pF_q$ [7, Entry 78 in \S
7.2.3], yet they 
are sufficiently general in nature and unify a number of known
results in 
literature [1, p.63], [27, pp.41-43], etc. Further the results
are obtained 
in a neat and compact form which may be used in investigating the
asymptotic 
behaviour of the generalized hypergeometric functions.
\par
From [8, p.101], we know that the generalized hypergeometric
function 
$_4F_3(\cdot)$ has the following Mellin-Barnes integral
representation
\begin{eqnarray}
&&\hspace*{-10mm}\Gamma\left[\begin{array}{l}a_1,a_2,a_3,a_4\\[2mm]
b_1,b_2,b_3\end{array}\right]
\
{}_4F_3\left(\begin{array}{l}a_1,a_2,a_3,a_4\\[2mm]b_1,b_2,b_3\end{array}
;z\right)\nonumber\\[2mm]
&&=\frac1{2\pi i}\int_{\Cs}\Gamma\left[\begin{array}{l}a_1+\xi,a_2+\xi,a_3+
\xi,a_4+\xi,-\xi\\[2mm]b_1+\xi,b_2+\xi,b_3+\xi\end{array}\right]
\ (-z)^\xi d\xi,
\end{eqnarray}
where the contour $\C$ is of Barnes type, extending from
$-i\infty$ to 
$+i\infty$ in the $\xi$-plane, curving if necessary, such that
the poles of 
gamma functions $\Gamma(a_1+\xi),\Gamma(a_2+\xi),\Gamma(a_3+\xi)$
and 
$\Gamma(a_4+\xi)$ lie to the left of the contour and the poles of

$\Gamma(-\xi)$ to its right.
\par
In what follows $m,n$ and $p$ are non-negative integers and the
notation
\begin{eqnarray}
&&\hspace*{-10mm}\Gamma\left[\begin{array}{l}a_1,a_2,a_3,a_4\\[2mm]
b_1,b_2,b_3\end{array}\right]=\frac{\Gamma(a_1)\Gamma(a_2)\Gamma(a_3)
\Gamma(a_4)}{\Gamma(b_1)\Gamma(b_2)\Gamma(b_3)}
\end{eqnarray}
is employed.
\par
The proofs of the results in \S 2 - \S 4 can be developed by an
appeal to the 
calculus of residues and following a procedure described in [5].
\\\\
\begin{flushleft}
{\bf 2. First Expansion}
\end{flushleft}
\setcounter{section}{2}
\par
We establish first the formula when two of the upper parameters
differ by 
integers:
\setcounter{equation}{0}
\begin{eqnarray}
&&\hspace*{-10mm}\Gamma\left[\begin{array}{l}a_1,a_1+m,a_3,a_4\\[2mm]
b_1,b_2,b_3\end{array}\right]\ {}_4F_3\left(\begin{array}{l}
a_1,a_1+m,a_3,a_4\\[2mm]b_1,b_2,b_3\end{array};z\right)\nonumber\\[2mm]
&&=(-z)^{-a_1}\sum^{m-1}_{\nu=0}\frac{z^{-\nu}}{\nu!}\Gamma\left[
\begin{array}
{l}m-\nu,a_3-a_1-\nu,a_4-a_1-\nu,a_1+\nu\\[2mm]b_1-a_1-\nu,b_2-a_1-\nu,
b_3-a_1-\nu\end{array}\right]\ z^{-\nu}\nonumber\\[2mm]
&&\hspace*{-15mm}+(-z)^{-a_1}\sum^\infty_{\nu=m}\frac{(-1)^m(-z)^{-\nu}}
{\nu!\thinspace(\nu-m)!}\Gamma\left[\begin{array}{l}a_1+\nu,a_3-a_1-\nu,
a_4-a_1-\nu\\[2mm]b_1-a_1-\nu,b_2-a_1-\nu,b_3-a_1-\nu\end{array}\right]
\ [A+\log(-z)],
\end{eqnarray}
where $m$ is a non-negative integer,
$a_1\ne0,-1,-2,\cdots;|\arg(-z)|<\pi;$
\begin{eqnarray}
&&\hspace*{-10mm}A=\psi(\nu+1)+\psi(\nu+1-m)+\psi(a_3-a_1-\nu)
+\psi(a_4-a_1-\nu)\nonumber\\[2mm]
&&-\psi(a_1+\nu)-\psi(b_1-a_1-\nu)-\psi(b_2-a_1-\nu)-\psi(b_3-a_1-\nu);
\end{eqnarray}
and $\psi(z)=d\{\log\Gamma(z)\}/dz.$
\\\par
For $a_4=b_3$ and $a_3=b_2$, the equation (2.1) gives rise to a
known result 
for the Gauss function $_2F_1$ [1, p.63].
\\\\
\begin{flushleft}
{\bf 3. Second Expansion}
\end{flushleft}
\setcounter{section}{3}
\setcounter{equation}{0}
\par
When three of the upper parameters differ by integers, there
holds the formula:
\begin{eqnarray}
&&\hspace{-5mm}\Gamma\left[\begin{array}{l}a_1,a_1+m,a_1+m+n,a_2\\[2mm]
b_1,b_2,b_3\end{array}\right]\
{}_4F_3\left(\begin{array}{l}a_1,a_1+m,
a_1+m+n,a_2\\[2mm]b_1,b_2,b_3\end{array};z\right)\nonumber\\[2mm]
&&=(-z)^{-a_1}\sum^{m-1}_{\nu=0}\frac{z^{-\nu}}{\nu!}\Gamma
\left[\begin{array}{l}m-\nu,m+n-\nu,a_1+\nu,a_2-a_1-\nu\\[2mm]
b_1-a_1-\nu,b_2-a_1-\nu,b_3-a_1-\nu\end{array}\right]\nonumber\\[2mm]
&&\hspace{5mm}+(-z)^{-a_1}(-1)^m\sum^{m+n-1}_{\nu=m}\frac{(-z)^{-\nu}}
{\nu!(\nu-m)!}\Gamma\left[\begin{array}{l}m+n-\nu,a_1+\nu,a_2-a_1-\nu\\[2mm]
b_1-a_1-\nu,b_2-a_1-\nu,b_3-a_1-\nu\end{array}\right]\nonumber\\[2mm]
&&\hspace{20mm}\cdot\ \{\theta+\log(-z)\}\nonumber\\[2mm]
&&\hspace{5mm}+(-z)^{-a_1}(-1)^n\sum^\infty_{\nu=m+n}\frac{z^\nu}
{\nu!(\nu-m)!(\nu-m-n)!}\nonumber\\[2mm]
&&\hspace{-13mm}\cdot\
\Gamma\left[\begin{array}{l}a_1+\nu,a_2-a_1-\nu\\[2mm]
b_1-a_1-\nu,b_2-a_1-\nu,b_3-a_1-\nu\end{array}\right]\
\left[B^2+B'+2\log(-z)
\thinspace B+\{\log(-z)\}^2\right],
\end{eqnarray}
where $m$ and $n$ are non-negative integers $a_1,\ a_1+m\ne0,-1,-2,\cdots$; 
and $|\arg(-z)|<\pi;$
\begin{eqnarray}
&&\hspace{-10mm}\theta=\psi(\nu+1)+\psi(\nu-m+1)+\psi(m+n-\nu)
-\psi(a_2-a_1-\nu)\nonumber\\[2mm]
&&-\psi(a_1+\nu)-\psi(b_1-a_1-\nu)-\psi(b_2-a_1-\nu)-\psi(b_3-a_1-\nu);\\[2mm]
&&\hspace{-10mm}B=\theta+\psi(\nu-m-n+1)-\psi(m+n-\nu);\\[2mm]
&&\hspace{-10mm}B'=3\psi'(1)+3\zeta(2,1)-\zeta(2,\nu+1)-\zeta(2,\nu-m+1)
-\zeta(2,\nu-m-n+1)\nonumber\\[2mm]
&&-\psi'(a_2-a_1-\nu)+\psi'(a_1+\nu)\nonumber\\[2mm]
&&-\psi'(b_1-a_1-\nu)-\psi'(b_2-a_1-\nu)-\psi'(b_3-a_1-\nu).
\end{eqnarray}
Here $\zeta(\cdot,\cdot)$ is the generalized Zeta function
defined by
\begin{eqnarray}
&&\zeta(s,\nu)=\sum^\infty_{n=0}(\nu+n)^{-s},\quad\nu\ne0,-1,-2,\cdots.
\end{eqnarray}
For a detailed account of generalized Zeta function
$\zeta(s,\nu)$, the 
reader is referred to the monograph by A. Erd\'elyi {\it et al}.
[1].
\\\\
\begin{flushleft}
{\bf 4. Third Expansion}
\end{flushleft}
\setcounter{section}{4}
\setcounter{equation}{0}
\par
When four of the upper parameters differ by integers, the
following result 
holds:
\begin{eqnarray}
&&\hspace{-5mm}\Gamma\left[\begin{array}{l}a_1,a_1+m,a_1+m+n,a_1+m+n+p\\[2mm]
b_1,b_2,b_3\end{array}\right]\nonumber\\[2mm]
&&\cdot\ {}_4F_3\left(\begin{array}{c}a_1,a_1+m,a_1+m+n,a_1+m+n+p\\[2mm]
b_1,b_2,b_3\end{array};z\right)\nonumber\\[2mm]
&&=(-z)^{-a_1}\sum^{m-1}_{\nu=0}\Gamma\left[\begin{array}{l}m-\nu,m+n-\nu,
m+n+p-\nu,a_1+\nu\\[2mm]b_1-a_1-\nu,b_2-a_1-\nu,b_3-a_1-\nu\end{array}\right]
z^{-\nu}\nonumber\\[2mm]
&&\hspace{5mm}+(-z)^{-a_1}\sum^{m+n-1}_{\nu=m}\frac{(-1)^mz^{-\nu}}
{\nu!(\nu-m)!}\Gamma\left[\begin{array}{l}m+n-\nu,m+n+p-\nu,a_1+\nu\\[2mm]
b_1-a_1-\nu,b_2-a_1-\nu,b_3-a_1-\nu\end{array}\right]\nonumber\\[2mm]
&&\hspace{20mm}\cdot\ \left[C+\log(-z)\right]\nonumber\\[2mm]
&&\hspace{5mm}+(-z)^{-a_1}\sum^{m+n+p-1}_{\nu=m+n}\frac{(-1)^nz^{-\nu}}
{2!\nu!(\nu-m)!(\nu-m-n)!}\nonumber\\[2mm]
&&\hspace{30mm}\cdot\
\Gamma\left[\begin{array}{l}m+n+p-\nu,a_1+\nu\\[2mm]
b_1-a_1-\nu,b_2-a_1-\nu,b_3-a_1-\nu\end{array}\right]\nonumber\\[2mm]
&&\hspace{40mm}\cdot\left[D'+D^2+2D\log(-z)+\left\{\log(-z)\right
\}^2\right]\nonumber\\[2mm]
&&\hspace{5mm}+(-z)^{-a_1}\sum^\infty_{\nu=m+n+p}\frac{(-1)^{m+p}(-z)^{-\nu}}
{\nu!(\nu-m)!(\nu-m-n)!(\nu-m-n-p)!}\nonumber\\[2mm]
&&\hspace{30mm}\cdot\ \Gamma\left[\begin{array}{l}a_1+\nu\\[2mm]
b_1-a_1-\nu,b_2-a_1-\nu,b_3-a_1-\nu\end{array}\right]\nonumber\\[2mm]
&&\hspace{3mm}\cdot\ \left[(E^3+3E'E+E'')+3\log(-z)(E'+E^2)
+3\{\log(-z)\}^2E+\{\log(-z)\}^3\right],
\end{eqnarray}
where $m,n$ and $p$ are non-negative integers,
$a_1,a_1+m,a_1+m+n\ne0,-1,-2,
\cdots,|\arg(-z)|<\pi$,
\begin{eqnarray}
&&\hspace{-10mm}C=\psi(\nu+1)+\psi(\nu-m+1)+\psi(m+n-\nu)+\psi(m+n+p-\nu)
-\psi(a_1+\nu)\nonumber\\[2mm]
&&-\psi(b_1-a_1-\nu)-\psi(b_2-a_1-\nu)-\psi(b_3-a_1-\nu);\\[2mm]
&&\hspace{-10mm}D=C+\psi(\nu-m-n+1)-\psi(m+n-\nu);\\[2mm]
&&\hspace{-10mm}D'=3\psi'(1)+3\zeta(2,1)-\zeta(2,\nu+1)-\zeta(2,\nu-m+1)
-\zeta(2,\nu-m-n+1)\nonumber\\[2mm]
&&-\psi'(m+n+p-\nu)+\psi'(a_1+\nu)\nonumber\\[2mm]
&&-\psi'(b_1-a_1-\nu)-\psi'(b_2-a_1-\nu)-\psi'(b_3-a_1-\nu);\\[2mm]
&&\hspace{-10mm}E=D+\psi(\nu-m-n-p+1)-\psi(m+n+p-\nu);\\[2mm]
&&\hspace{-10mm}E'=D'+\psi'(1)+\zeta(2,1)+\psi'(m+n+p-\nu)
-\zeta(2,\nu-m-n-p+1);\\[2mm]
&&\hspace{-10mm}E''=2\left[2\psi''(1)+4\zeta(3,1)-\zeta(3,\nu+1)
-\zeta(3,\nu-m+1)\right.\nonumber\\[2mm]
&&\hspace{10mm}\left.-\zeta(3,\nu-m-n+1)-\zeta(3,\nu-m-n-p+1)\right]
\nonumber\\[2mm]
&&-\psi''(a_1+\nu)-\psi''(b_1-a_1-\nu)-\psi''(b_2-a_1-\nu)-\psi''
(b_3-a_1-\nu).
\end{eqnarray}
%
%\\\par
%
\newpage
{\bf Remark 1.} \ It may be mentioned here that the results given
earlier by 
Saxena and Kalla [26, pp.41-43] follow from (2.1) and (3.1).
\\\par
{\bf Remark 2.} \ Though the formulae (2.1), (3.1) and (4.1) are
special 
cases of a general one for $_pF_q$ listed in [7, \S 7.2.3.78],
they are still 
sufficiently general and unify known results in [1, p.63], [27,
pp.41-43] 
etc.
\\\par
{\bf Acknowledgement.} \ The present investigation was supported,
in part, by 
the Ministry of 
Education of Japan under Grant 09640238 and the University Grants
Commission of India. Thanks are due to 
Professor H.M. Srivastava for giving some useful suggestions.
\vspace{1cm}
\begin{center}{\bf References}\end{center}
\small
\begin{enumerate}
\item[{[1]\ }]A. Erd\'elyi, W. Magnus, F. Oberhettiger and F.G.
Tricomi: 
{\it Higher Transcendental Functions}, Vol.1, McGraw-Hill, New
York, 1953.
\vspace{-2mm}
\item[{[2]\ }]W. Magnus, F. Oberhettiger and R.P Soni: {\it
Formulas and 
Theorems for the Special Functions of Mathematical Physics},
Third Enlarged 
Ed., Springer-Verlag, Berlin-Heidelberg, 1966.
\vspace{-2mm}
\item[{[3]\ }]A.M. Mathai: A few remarks on the exact
distributions of 
certain multivariate statistics - II. {\it Multivariate
Statistical 
Inference}, pp. 169-181, Kabe, D.G. and Gupta, R.P. (Eds.), North
Holland, 
Amsterdam, 1972.
\vspace{-2mm}
\item[{[4]\ }]A.M. Mathai and H.J. Haubold: {\it Modern Problems
in Nuclear 
and Neutrino Astrophysics}, Akademie-Verlag, Berlin, 1988.
\vspace{-2mm}
\item[{[5]\ }]A.M. Mathai and R.K. Saxena: {\it Generalized
Hypergeometric 
Function with Applications in Statistics and Physical Sciences},
Lecture Notes 
Math., Vol. 348, Berlin, 1973.
\vspace{-2mm}
\item[{[6]\ }]A.M. Mathai and R.K. Saxena: {\it The $H$-function
with 
Applications in Statistics and Other Disciplines}, Halsted Press,
New 
York-London-Sydney-Toronto, 1978.
\vspace{-2mm}
\item[{[7]\ }]A.P. Prudnikov, Yu.A. Brychkov and O.I. Marichev:
{\it Integral 
and Series}, Vol.3, {\it More Special Functions}, Gordon and
Breach, New 
York-London, 1989.
\vspace{-2mm}
\item[{[8]\ }]E.D. Rainville: {\it Special Functions}, Chelsea
Pub., New York, 
1960.
\vspace{-2mm}
\item[{[9]\ }]M. Saigo: On a property of the Appell
hypergeometric function 
$F_1$, {\it Math. Rep. Kyushu Univ.} {\bf 12}(1980), 63-67.
\vspace{-2mm}
\item[{[10]\ }]M. Saigo: On properties of the Appell
hypergeometric functions 
$F_2$ and $F_3$ and the generalized Gauss function $_3F_2$, {\it
Bull. Centr. 
Res. Inst.}, {\it Fukuoka Univ.} {\bf 66}(1983), 27-32.
\vspace{-2mm}
\item[{[11]\ }]M. Saigo: On properties of the Lauricella
hypergeometric 
function $F_D$, {\it Bull. Centr. Res. Inst.}, {\it Fukuoka
Univ.}
{\bf 104}(1988), 13-31.
\vspace{-2mm}
\item[{[12]\ }]M. Saigo: On properties of hypergeometric
functions of three 
variables, $F_M$ and $F_G$, {\it Rend. Circ. Mat. Palermo}, {\it
Serie II} 
{\bf 37}(1988), 449-468.
\vspace{-2mm}
\item[{[13]\ }]M. Saigo: On behaviors of hypergeometric series of
one, two, 
three and $n$-variables near boundaries of their convergence
regions of 
logarithmic case, (Kyoto, 1989), {\it S\^urikaiseki-kenky\^usho
K\^oky\^uroku} 
{\bf 714}(1990), 91-109.
\vspace{-2mm}
\item[{[14]\ }]M. Saigo: The asymptotic behaviors of triple
hypergeometric 
series $F_G,F_K,F_N$ and $F_R$ near boundaries of their
convergence regions 
(Kyoto, 1989), {\it S\^urikaisekikenky\^usho K\^oky\^uroku} {\bf
714}(1990), 
110-126.
\vspace{-2mm}
\item[{[15]\ }]M. Saigo: Asymptotic behaviors of the Horn
hypergeometric 
series $G_1$, $G_2$, $H_1$, $H_2$, $H_4$, $H_6$, $H_7$ near
boundaries of 
their convergence regions, {\it Integral Transform. Spec. Funct.}

{\bf 4}(1996), 249-262.
\vspace{-2mm}
\item[{[16]\ }]M. Saigo, O.I. Marichev and Ngyen Thanh Hai:
Asymptotic 
representations of Gaussian series $_2F_1$, Clausenian series
$_3F_2$ and 
Appell series $F_2$ and $F_3$ near boundaries of their
convergence regions, 
{\it Fukuoka Univ. Sci. Rep.} {\bf 19}(1989), 83-90.
\vspace{-2mm}
\item[{[17]\ }]M. Saigo and Ngyen Thanh Hai: The asymptotic
representations 
of Appell series $F_1$ and Horn series $H_3$ near boundaries of
their 
convergence regions, {\it Fukuoka Univ. Sci. Rep.} {\bf
21}(1991), 17-23.
\vspace{-2mm}
\item[{[18]\ }]M. Saigo and Ngyen Thanh Hai: Asymptotic
representations of 
the Horn series $G_1$ and $G_2$ near boundaries of their
convergence regions, 
{\it Fukuoka Univ. Sci. Rep.} {\bf 22}(1992), 1-9.
\vspace{-2mm}
\item[{[19]\ }]M. Saigo and H.M. Srivastava: The behaviors of the
Appell 
double hypergeometric series $F_4$ and certain Lauricella triple 
hypergeometric series near the boundaries of their convergence
regions, 
{\it Fukuoka Univ. Sci. Rep.} {\bf 19}(1989), 1-10.
\vspace{-2mm}
\item[{[20]\ }]M. Saigo and H.M. Srivastava: The asymptotic
behaviour of 
Lauricella's second hypergeometric series $F_B^{(n)}$ in $n$
variables near 
the boundary of its convergence region, {\it Fukuoka Univ. Sci.
Rep.} 
{\bf 19}(1989), 71-82.
\vspace{-2mm}
\item[{[21]\ }]M. Saigo and H.M. Srivastava: The behaviors of the

zero-balanced hypergeometric series $_pF_{p-1}$ near the boundary
of its 
convergence region, {\it Proc. Amer. Math. Soc.} {\bf 110}(1990),
71-76.
\vspace{-2mm}
\item[{[22]\ }]M. Saigo and H.M. Srivastava: Some asymptotic
formulas 
exhibiting the behaviors of the triple hypergeometric series
$F_S$ and $F_T$ 
near the boundaries of their convergence regions, {\it Funk.
Ekvac.} 
{\bf 34}(1991), 423-448.
\vspace{-2mm}
\item[{[23]\ }]M. Saigo and H.M. Srivastava: A theorem on the
asymptotic 
behavior of Lauricella's first hypergeometric series $F_A^{(n)}$
near the 
boundary of its convergence region, {\it Complex Analysis and
Generalized 
Functions}, Varna, September 15-21 (1991), 258-267.
\vspace{-2mm}
\item[{[24]\ }]M. Saigo and H.M. Srivastava: The asymptotic
behaviour of the 
Srivastava hypergeometric series $H_C$ near the boundary of its
convergence 
region, {\it Kyungpook Math. J.} {\bf 32}(1992), 153-158.
\vspace{-2mm}
\item[{[25]\ }]M. Saigo and H.M. Srivastava: The behaviour of the
fourth type 
of Lauricella's hypergeometric series in $n$ variables near the
boundaries of 
its convergence region, {\it J. Austral. Math. Soc.$,$ Ser. A}
{\bf57}(1994), 
281-294.
\vspace{-2mm}
\item[{[26]\ }]R.K. Saxena and S.L. Kalla: Expansion of a $_3F_2$
when the 
upper parameters differ by an integer, {\it Rev. Tec. Ing. Univ.
Zulia} 
{\bf 9}(1986), 41-43.
\end{enumerate}
\end{document}